\documentclass[a4,12pt]{article}

\usepackage{xcolor} % Required for specifying colors by name
\definecolor{ocre}{RGB}{52,177,201} % Define the orange color used for highlighting throughout the book

\usepackage{diagbox}
\usepackage{ stmaryrd }

\usepackage{listings}[language=Bash]
\usepackage{afterpage}

\usepackage{avant} % Use the Avantgarde font for headings
\usepackage{mathptmx} % Use the Adobe Times Roman as the default text font together with math symbols from the Sym­bol, Chancery and Com­puter Modern fonts

\usepackage{microtype} % Slightly tweak font spacing for aesthetics
\usepackage[utf8]{inputenc} % Required for including letters with accents
\usepackage[T1]{fontenc} % Use 8-bit encoding that has 256 glyphs

\usepackage{float}
\usepackage[caption = false]{subfig}
\usepackage[final]{graphicx}

\usepackage{mdwlist}
\usepackage{multirow}
\usepackage{siunitx}

\graphicspath{{img/}} % Specifies the directory where pictures are stored

\usepackage{subfig}

\usepackage[english]{babel} % English language/hyphenation

\usepackage{enumitem} % Customize lists
\setlist{nolistsep} % Reduce spacing between bullet points and numbered lists

\usepackage{booktabs} % Required for nicer horizontal rules in tables

\definecolor{Navy}{RGB}{000,000,128}
\definecolor{DarkRed}{RGB}{139,000,000}
\definecolor{DarkOrange}{RGB}{205,102,000}
\definecolor{DarkGreen}{RGB}{000,100,000}
\definecolor{Beige}{RGB}{250,250,240}
\definecolor{applegreen}{rgb}{0.55, 0.71, 0.0}
\definecolor{palecopper}{rgb}{0.85, 0.54, 0.4}
\definecolor{platinum}{rgb}{0.9, 0.89, 0.89}
\definecolor{CAUpurple}{RGB}{153,000,153}
\definecolor{CAUblue}{RGB}{051,051,153}
\definecolor{DarkPurple}{rgb}{0.07, 0.04, 0.56}

\definecolor{codegreen}{rgb}{0,0.6,0}
\definecolor{codegray}{rgb}{0.5,0.5,0.5}
\definecolor{codepurple}{rgb}{0.58,0,0.82}
\definecolor{backcolour}{rgb}{0.95,0.95,0.92}

\lstdefinestyle{mystyle}{
	backgroundcolor=\color{backcolour},   
	commentstyle=\color{codegreen},
	keywordstyle=\color{magenta},
	numberstyle=\tiny\color{codegray},
	stringstyle=\color{codepurple},
	basicstyle=\ttfamily\footnotesize,
	breakatwhitespace=false,         
	breaklines=true,                 
	captionpos=b,                    
	keepspaces=true,                 
	numbers=left,                    
	numbersep=5pt,                  
	showspaces=false,                
	showstringspaces=false,
	showtabs=false,                  
	tabsize=2
}

\usepackage{fancyhdr} % Required for header and footer configuration

\usepackage{amsmath,amsfonts,amssymb,amsthm} % For math equations, theorems, symbols, etc

\usepackage{hyperref}

\hypersetup{
    bookmarks=true,
    unicode=false,
    pdftoolbar=true,
    pdfmenubar=true,
    pdffitwindow=false,
    pdfstartview={FitH},
    pdftitle={My title},
    pdfauthor={B. Philippi},
    pdfsubject={Subject},
    pdfcreator={Creator},
    pdfproducer={Producer},
    pdfkeywords={keyword1} {key2} {key3},
    pdfnewwindow=true,
    colorlinks=false,
    linkbordercolor={0 0 1},
    linkcolor=blue,
    citecolor=green,
    filecolor=magenta,
    urlcolor=cyan
}

\usepackage[left=3cm,right=3cm,top=2cm,bottom=2cm,includeheadfoot]{geometry}

\title{The Parareal Algorithm Applied to the FESOM 2 Ocean Circulation Model}
\author{B. $\text{Philippi}^1$, T. $\text{Slawig}^2$} 

\begin{document}

\maketitle

\begin{abstract}
	In this work the parallel-in-time algorithm Parareal was applied to the ocean-circulation and sea-ice model FESOM2 developed by the Alfred-Wegener Institut (AWI). The climate model provides one time integration method and hence, the coarse and fine propagators were defined by time step width. The coarse method was executed at the CFL condition limit, while fine step-sizes where gradually refined. As a first assessment of the performance of Parareal a low-resolution test mesh was used with default settings provided by the AWI. An introduction to FESOM2 and the straightforward implementation of Parareal on the DKRZ cluster is given. The evaluation of numerical results for different simulation intervals and fine propagator configurations shows strong dependence on the simulation time and time step-size of the fine propagator. Increasing the latter leads to stagnation and eventually divergence of the Parareal algorithm.
	
\end{abstract}

\vspace{0.25cm}
\small
\textbf{${}^1$ \textit{Christian-Albrecht-Universität Kiel, Dept. of Computer Science, bph@informatik.uni-kiel.de}}

\textbf{${}^2$ \textit{Christian-Albrecht-Universität Kiel, Dept. of Computer Science, ts@informatik.uni-kiel.de}}
\normalsize

\newpage

\section{Introduction}

With the growing development of modern supercomputers' processing power, the demand for efficient algorithms for parallelization in space rises accordingly. However, not all computational applications do benefit from increased spatial resolution, and hence the vast capacities of high-performance clusters (HPC) are far from exploited. Examples for such applications can be found in geophysical models, where large scale physics for atmosphere, ocean circulation and sea-ice evolution over long time periods are simulated. Due to the high complexity in geophysical phenomenons, the models are focused on a specific physical domain and optimized with respect to the chosen spatial resolution. This often results in working configurations that are then used for hind- and forecast simulations or for coupling to models addressing other physical aspects of the earth's physics. Due to parameterized design of the models, time-demanding calibration processes are necessary for respective spatial resolutions. Thus, runtime reduction by domain decomposition techniques will necessarily saturate at some point. The concept of parallel-in-time algorithms extends beyond classical parallelization approaches by decomposing the given problem along the time axis. As one of these approaches, the \textit{Parareal} algorithm is well-known and has been extensively studied since it was first introduced by \cite{Maday2001}. Although we restrict ourselves to Parareal in this study, other notable contributions to the parallel-in-time community have been made, e.g.,  by the PFASST \cite{PFASST} and PITA \cite{FarhatEtAl2003} approaches. For a broad overview of time-parallel time integration methods and their history the reader is referred to \cite{Gander2015Review}.

Although the algorithm was initially introduced for ordinary differential equations (ODEs), the field of applications has since been extended towards linear and non-linear partial differential equations (PDEs). Fluid dynamic simulations described by the Navier-Stokes equations still represent a challenge for time-parallel algorithms, especially when the flow enters the turbulent regime. For the classical driven cavity flow test case with additional obstacles placed within the domain, a successful implementation of Parareal was presented in \cite{Ruprecht2014}. A vortex shedding flow with a Reynolds number of $Re=1000$ was assessed in an early study in \cite{FarhatEtAl2003}, where the authors observed that Parareal's stability highly depended on the spatial dimension and the used time integration methods. Further, it was confirmed that with increasing Reynolds number and spatial resolution the convergence rate of Parareal would decline in cases of the fully solved Navier-Stokes equations \cite{Ruprecht2015}. For direct numerical simulation (DNS, i.e., the computation of the full Navier-Stokes equations without any turbulence modeling) with Reynolds numbers higher than 10000, the approximation of states at specific points in time by Parareal fails. However, it was still possible to approximate main statistical quantities of flows in turbulent regimes, as it has been done for the decay of homogeneous turbulence in a three-dimensional test case in \cite{Lunet2018}. 

In the context of numerical ocean science, a time-parallel execution of Navier-Stokes-based ocean models can be of interest. As the evolution of geophysical phenomenons is computationally demanding, the question for runtime reduction arises. Characteristic diagnostic variables, e.g., the sea surface temperature (SST) and mean ocean circulation (MOC), are important quantities that are derived from local states in time. To the authors' knowledge, the application of Parareal to an ocean simulation software has been investigated once by \cite{Liu2008}, where the Princeton ocean model was used. The test case concerned the Bohai Sea representing a marginal part of the world's oceans. In contrast, we applied the Parareal algorithm to the sea-ice ocean circulation model FESOM2 as a first assessment in this work. The ocean model FESOM2 was developed at the Alfred-Wegener-Institut, Helmholtz Centre for Polar and Marine Research (AWI), and simulates the long-time evolution of sea-ice and ocean dynamics \cite{FESOM2oce}. The set of main equations represent the Navier-Stokes equations with turbulence models adjusted to the respective mesh resolution \cite{FESOM2standardSim}. Since sea-ice dynamics represent a driving force in the ocean circulation, FESOM2 contains the implementation of a simplified ice model \cite{FESOM2ice}. The focus of this work was set to analyze Parareal's convergence when it comes to the approximation of states at specified points in time. We chose a low resolution test mesh in accordance with the findings of the convergence study in \cite{Ruprecht2015} which confirmed good performance of Parareal for low spatial resolution. Nevertheless, we expected Parareal to struggle with stability issues since the ocean dynamics represents an advection-dominant problem, compare \cite{GanderVandewalle},\cite{Ruprecht2015}. In this paper we will introduce a straightforward implementation of Parareal for FESOM2 on the HPC of the \textit{Deutsches Klimarechenzentrum} (DKRZ). Afterwards we will discuss the results for a default test setting on the so-called \textit{PI} mesh with low spatial resolution. The conducted numerical experiments presented in this paper serve the purpose of assessing the performance of the time-parallel algorithm for an established ocean model for climate research. With the completion of this study, the applicability can be estimated and next steps for an efficient implementation of Parareal to standard configurations of the ocean modeling community may be considered.

\section{Parareal Algorithm}

Parareal is a parallel-in-time algorithm for the approximation of initial value problems (IVPs) of the form:
\begin{equation}
u_t(t) \; = \; f(u(t),t), \quad t \in [0,T], \quad u(0) \; = \; u_0.
\label{EQ:IVP}
\end{equation}
The algorithm iterates approximations to the solution of \eqref{EQ:IVP} by applying two time integration methods that  differ in numerical accuracy and computational cost. The accurate solver is commonly denoted as the \textit{fine propagator} $F$ and the fast solver as the \textit{coarse propagator} $G$. The temporal domain $[0,T]$ is split into $N_t$ sub-intervals $\Delta T = [t_{n},t_{n+1}]$, $n=0,\dots,N_t-1$. To simplify the notation and describe the basic structure of the algorithm, we assume that one application of each of the two propagators  approximates  a solution from a given state $U_{n}\approx u(t_{n})$ at time $t_{n}$ over the time slice with length $\Delta T$ to time $t_{n+1}$, i.e.:
\begin{equation}
\begin{split}
F_{n+1}\;&:=\; F(U_{n}) \approx u(t_{n+1}), \\
G_{n+1}\;&:=\; G(U_{n}) \approx u(t_{n+1}).
\end{split}
\end{equation}
Here, any direct dependency of $F$ and $G$ on the time instant $t_{n}$ is omitted in the notation.
In reality, the time integration of both propagators over the time slice $[t_{n},t_{n+1}]$ may consist of  several time steps of the used numerical time integrator.
The classical approach is solving the given problem over all time slices in a sequential way by applying $F$ from $t_0=0$ to $t_{N_t}=T$. In contrast, the concept of the parallel-in-time algorithm Parareal is to iterate by solving IVPs on the sub-intervals in parallel.
For this purpose, the algorithm starts by a sequential run with the coarse propagator $G$ over the whole time interval, giving a first iterate (the superscript representing the iteration counter):
\begin{equation}
\label{EQ:coarseinitialrun}
\begin{split}
U^{0}_{0}\;&:= \;u_{0},\\
U^{0}_{n+1} \; &:= \; G^{0}_{n+1}\;:=\;G(U^{0}_{n}) , \quad n=0,\ldots, N_{t}-1.
\end{split}
\end{equation}
Providing these values on the matching points of the time slices by the fast solver $G$ allows for a time-parallel computation of the time-consuming terms involving the fine solver.
Given the $k$-th iterate, one  iteration  then contains a parallel computation of the values 
\begin{equation}
\label{EQ:FG}
 F^{k}_{n+1}\;:=\;F(U^k_{n}), \quad  n=0,\ldots, N_{t}-1,
\end{equation}
with $U^{k}_{0}=u_{0}$ for all $k$. 
Afterwards, a sequential correction run with the coarse propagator is executed over the whole time interval, where the difference of the fine and coarse results of the last iteration are added  as corrections in every time step. This gives the next iterate:
\begin{equation}
\begin{split}
U^{k+1}_{0}\;&:= \;u_{0},\\
G^{k+1}_{n+1}\;&:=\;G(U^{k+1}_{n}),\\
U^{k+1}_{n+1} \; &:= \; G(U^{k+1}_{n}) \; + \; F^k_{n+1} \; - \; G^k_{n+1}, \quad n=0,\ldots, N_{t}-1.
\end{split}
\end{equation}
 By applying the computationally expensive method $F$ parallel-in-time a speed-up is generated if $k \ll N_t$. It can be shown that with $k \rightarrow N_t$, the iterate $U^k_n$ will converge towards the fine solution. For further interest in a convergence study and interpretation of the algorithm the reader is referred to \cite{GanderVandewalle}.

\subsection{Speed-up estimate}

Parareal allows for speed-up by enabling a further level of parallelization. Investing more computational resources into the given problem makes sense, when runtime reduction by domain decomposition techniques is saturated. In order to justify the additional amount of resources an a-priori speed-up estimate of the algorithm is useful. Let $\tau_G$ and $\tau_F$ denote the runtimes necessary to evaluate $G$ and $F$ over one sub-interval, respectively. We assume that the runtime for each time interval  of the fine propagator $\tau_F$ can be expressed in dependency of $\tau_G$, such that $\tau_F = m \cdot \tau_G$, with $m$ as a model factor that can be estimated or derived from runtime tests. The  serial computation runtime with the fine propagator is given as
\begin{equation*}
R_S(N_t,\tau_F)\;=\;N_t \cdot \tau_F
\;=\;N_t \cdot m \cdot \tau_G,
\end{equation*}
 and the Parareal runtime as
 \begin{equation*}
 \begin{split}
 R_P(N_t,k,\tau_F,\tau_G)
 \;=\;
 N_t \cdot \tau_G \; +\; k \cdot N_t \cdot \tau_G \; +\; k \cdot \tau_F
 \;=\;
(k+1) \cdot N_t \cdot \tau_G \; +\; k \cdot m\cdot \tau_G,
 \end{split}
 \end{equation*}
 both over the whole time interval $[0,T]$.
 The expected speedup is denoted as the ratio
\begin{equation*}
S(k,N_t,m) \; =  \; \frac{R_S(N_t,\tau_F)}{R_P(N_t,k,\tau_F,\tau_G)}
\;= \; \frac{N_t \cdot m}{(k+1)\cdot N_t  \; +\; k m} 
\;=\;  \frac{m}{k  + 1 \; + \; \frac{k m}{N_t} }  
\;=\; \frac{1}{\frac{k+1}{m} \; + \; \frac{k}{N_t}}.
\end{equation*}
Thus, a simplified, rough theoretical upper bound for $S$ is given by: 
\begin{equation}
	S(k,N_t,m) \; \leq \; \min \left\{ \frac{m}{k+1} \; , \; \frac{N_t}{k}  \right\} \; .
	\label{EQ:speedup}
\end{equation}
Basically it says: The faster the coarse propagator is compared to the fine one (factor $m$) and
the higher the number of available parallel processes (and thus time slices $N_{t}$) is compared to the needed iterations, the better Parareal works.  
The speed-up estimate given here does not account for data file processing and node communication. These quantities can be evaluated with monitoring software while executing the algorithm on a cluster. Consequently  Eq. \eqref{EQ:speedup} does provide a speed-up measure overestimating the true runtime reduction on a HPC environment. Nevertheless, as an indicator for possible speedups the derived a-priori measure is considered sufficient.  

\section{FESOM 2}

In this study, version 2.0 of the Finite-volumE Sea ice-Ocean circulation Model (FESOM2) is chosen to evaluate the parallelization-in-time concept of Parareal. It is formulated on unstructured finite volume meshes. The mayor advantage over its Finite-Element-Method based predecessor FESOM 1.4 is a computational efficiency comparable to ocean models based on structured meshes. The ocean mesh is split in horizontal planes that are aligned in vertical direction throughout the ocean. By doing so, the numerical efficiency and scalability for parallel computation on HPCs was significantly improved, see \cite{FESOM2performance}. 

FESOM2 provides linear and non-linear free-surface models, turbulent eddy modeling by the Gent-McWilliams parameterization and an isoneutral Redi diffusion. For a detailed description of the solver and the implementation of all its numerical models, see \cite{FESOM2oce}, \cite{FESOM2standardSim}. All variables of ocean and sea ice are intertwined with each other and need to be investigated as an entity to understand the complex physical behavior of the system. However, a evaluation of Parareal with respect to all variables would go beyond the scope of this study. And hence, during the course of the numerical experiments the focus is set on the investigation of zonal velocity $u$, temperature $T$ and salinity $S$. We figured the feasibility of Parareal for FESOM2 as a first assessment is best estimated by these variables, since they are computed by the core set of equations. Evaluation of the ocean's surface temperature, salinity distributions and the meridional overturning circulation derived from the velocity components are important diagnostic variables in maritime geoscience. 

By splitting the mesh in vertical aligned layers, the horizontal and vertical components of the governing equations are treated separately and the set of equations obtained is denoted accordingly. The horizontal momentum conservation is given by:
\begin{equation}
	\partial_t \mathbf{u} \; + \; f \mathbf{e}_z \times \mathbf{u} \; + \; \left( \mathbf{u} \cdot \nabla_h \; + \; w \partial_z \right) \mathbf{u} \; + \; \nabla_h \frac{p}{\rho_0} \; = \; D_h \mathbf{u} \; + \; \partial_z \nu_z \partial_z \mathbf{u} \; ,
\end{equation}
where $\mathbf{u} = (u,v)$ denotes the horizontal velocity, $f$ the Coriolis parameter, $p$ the pressure, $\nu_z$ the vertical viscosity,  and $\rho_0$ the reference density. The operators $\nabla_h = (\partial_x,\partial_y)$ and $D_h$ address the horizontal gradient and the horizontal viscosity term. Terms concerning the vertical contributions carry the sub-index $z$. The vertical velocity is denoted as $w$. The tracer equations for temperature $T$ and salinity $S$ are given by:
\begin{equation}
	\begin{split}
		\partial_t T \; + \; \nabla_h \cdot (\mathbf{u} T) \; + \; \partial_z (w T) \; &= \; \nabla_h \cdot \mathbf{K} \nabla_h T \; , \\
		\partial_t S \; + \; \nabla_h \cdot (\mathbf{u} S) \; + \; \partial_z (w S) \; &= \; \nabla_h \cdot \mathbf{K} \nabla_h S \; ,
	\end{split}
	\label{EQ:tracer}
\end{equation}
with the diffusivity tensor $\nu_v$ for mixing in the ocean. FESOM2 solves the main equations on a spherical finite-volume mesh and therefore is dependent on parameterizations to account for the complexity of the general ocean dynamics. The large scale sea-ice dynamics are described by the sub-model FESIM, the standard sea-ice model coming with FESOM2, which is  a zero-layer thermodynamical model, see \cite{FESOM2ice}. A detailed introduction to all parts of FESOM2 lies beyond the scope of this paper. The reader is referred to the given references for the implementation and evaluation of said models.

FESOM2 saves simulation data after a successful finalization in restart files for sea ice and ocean related variables and, if desired, selected post-processing quantities in separate files. All simulation data are gathered within netCDF files. The Network Common Data Format (\textit{netCDF}, \cite{netcdf}) is a self-describing data format widely used in climatology, meteorology and oceanography. For the netCDF format, software to store, access and manipulate arrays of scientific data sets is available. Additionally to the restart files, a simulation-time file with respect to the latest time stamp of the simulation is provided: \\
\begin{lstlisting}
	output/ 
		fesom.1948.oce.restart.nc
		fesom.1948.ice.restart.nc
		fesom.clock
\end{lstlisting} 
\vspace{.4cm}

These files contain all ocean and sea ice states of the last time step and therefore are of special interest for the implementation of Parareal. As the name assignment suggests, FESOM2 allows for continuation the simulation from the restart files to propagate further in time. For the sake of clarity, we restrict ourselves to the restart files, since just they will be used for the parallel-in-time algorithm. 

The various simulation parameters for the ocean and sea ice models can be defined in several configuration files. For the purpose of this study, a default setting from the AWI was taken. It is obtainable on \hyperlink{fesom2}{https://github.com/FESOM/fesom2}. Alongside the default settings, we obtained mesh files and atmospheric forcing files from the FESOM2 github repository. In addition to the default settings, an atmospheric forcing has to be chosen. In this way the coupled simulation of ocean and atmosphere can be avoided. Instead,   a data set that describes the atmospheric influence on ocean and sea-ice is provided. Forcing data sets contain information for air temperature, solar heat, precipitation and wind velocities above the ocean. Ocean simulation runs with such forcing data  are referred to as hindcast experiments, since they rely on data sets that are composed over years of real-life measurements. For all experiments the CORE-II forcing was used, which covers the 60-year time period from 1948 to 2007, see \cite{core2}.

During this study the coarsest mesh available, the PI mesh, was chosen to evaluate the performance of Parareal. The mesh contains 3140 nodes for each of the 48 horizontal layers. The mesh files are freely obtainable alongside the aforementioned default parameter settings on github (\url{https://github.com/FESOM}). The paths to the locations on the computing system of mesh and forcing files as well as the FESOM2 executable are given among the parameter settings within several \textit{namelist} files. The simulation folder containing all information and the output directory have to be created beforehand: \\
\begin{lstlisting}
	run_folder/
		fesom.x -> $PATH_TO_EXECUTABLE/fesom.x
		namelist.config
		namelist.forcing
		namelist.ice
		namelist.io
		namelist.oce
		output/
\end{lstlisting} 
In the file \texttt{namelist.config} settings for simulation time, time step-size, the location of the mesh file and the amount of processors for domain decomposition are defined. Therefore this file is of particular interest for the implementation of the Parareal algorithm. The file \texttt{namelist.forcing}  contains links to the forcing data set, \texttt{namelist.io} controls the output frequency of diagnostic variables, \texttt{namelist.oce} and \texttt{namelist.ice} allow for manipulation of the ocean and sea-ice parameterizations.

\begin{figure}[H]
	\centering
	\includegraphics[width=.45\textwidth]{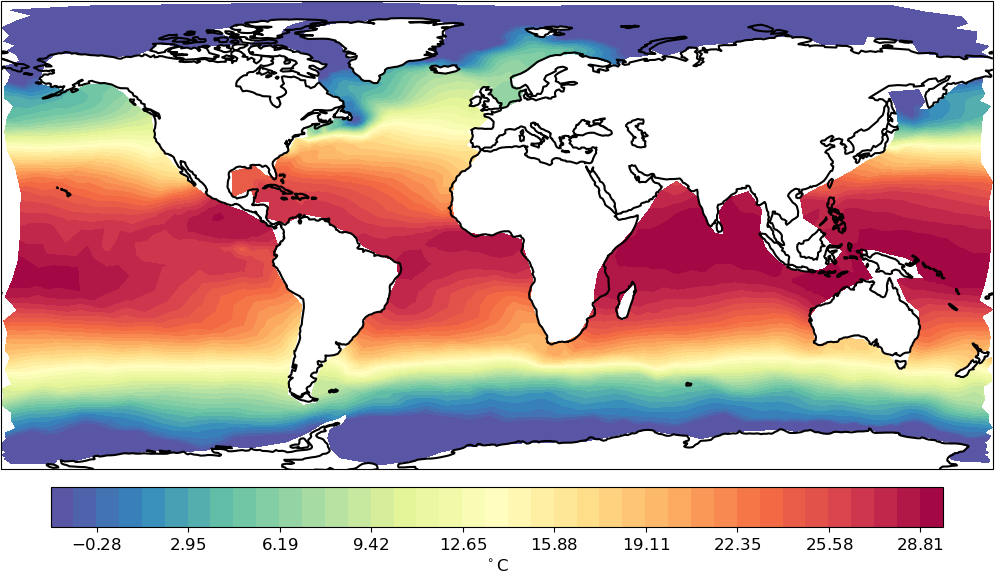} \quad
	\includegraphics[scale=.2]{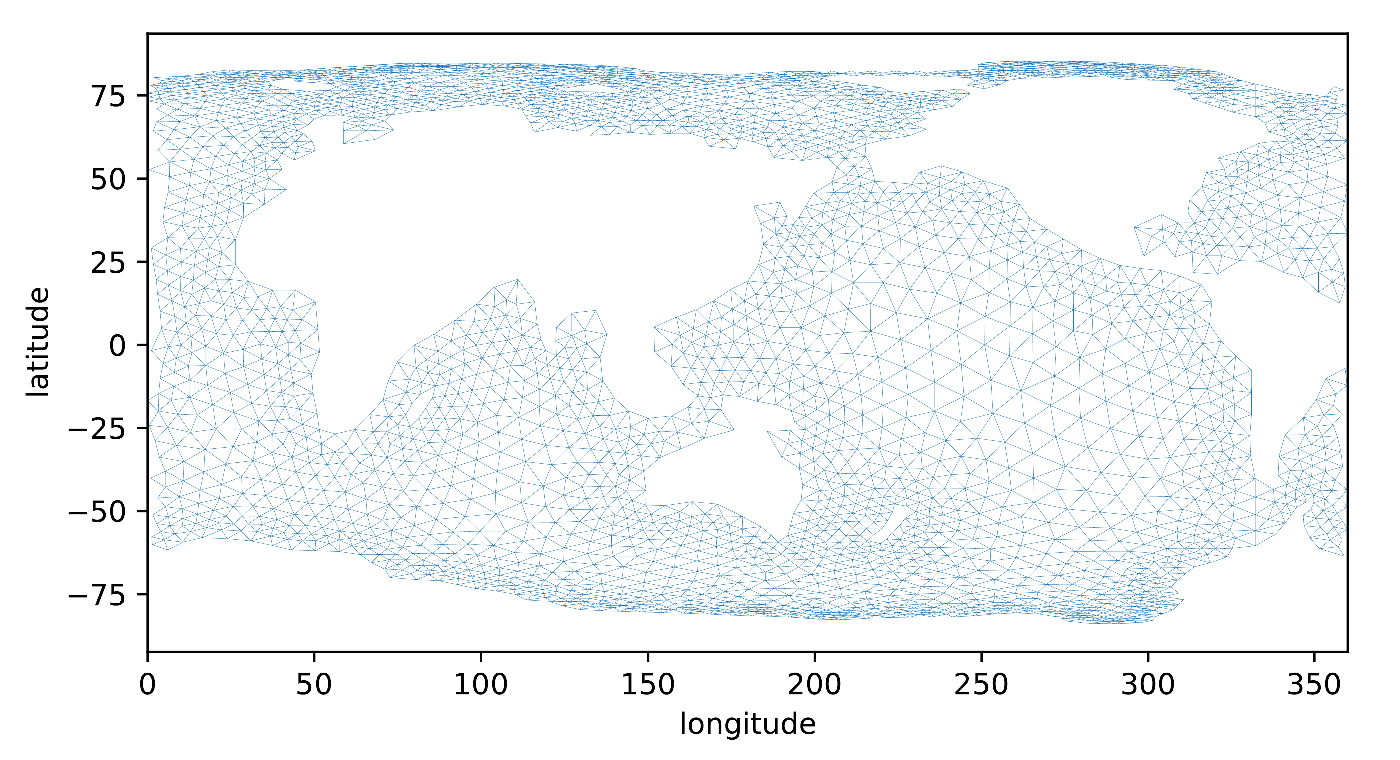}
	\caption{Left: Example temperature distribution computed on the PI mesh with CORE2 forcing. Right: Visualization of the finite-volume discretization.}
	\label{FIG:TempMeshEx}
\end{figure}

In Fig.\ref{FIG:TempMeshEx} an example temperature distribution computed by FESOM2 and the underlying finite-volume PI mesh are shown. Due to its low spatial resolution, some areas could not be modelled, for example the Baltic and Mediterranian Sea. Small passages, like the Strait of Gibraltar, pose numerical problems when one tries to incorporate them into coarse meshes.

\vspace{1cm}

\section{Parareal implementation for the DKRZ Cluster}

The DKRZ (Deutsches Klimarechenzentrum) provides the petascale High Performance Computing system Mistral. Cluster management and job scheduling system are operated with \textit{Slurm}, an open source workload manager. For each job a maximum of 512 nodes with 24 processors each can be requested. The instructions for each job are gathered within a batch script which is submitted to the workload manager. The so-called \textit{sbatch} script contains all necessary steps and settings for Parareal to be executed. 

In order to submit a job, the computing partition, the amount of nodes, number of tasks per node and the estimated runtime need to be defined preceding the executable commands within the script. If 12 tasks are to be executed in parallel with 24 processors each, the total amount of resources must be defined in the sbatch options at the very beginning of the script. The workload manager will process the sbatch script in a serial way. However, for Parareal it is necessary to execute parallel sub-steps with a subset of requested resources. The Slurm command \texttt{srun} allows for such substeps in any desired configuration, and therefore provides time-parallel execution of the fine propagator. Furthermore, the batch programming language contains control structures like \texttt{if, for, while}. Hence, scripting the algorithm's steps comes in a straight-forward manner. A pseudo-code representation of Parareal is given in Listing \ref{lst:Parareal}. \\ 

\begin{lstlisting}[language=octave,caption={Pseudo-Code of the Parareal algorithm},label={lst:Parareal}]

	# initialization run of the coarse propagator in serial

	for n=0:Nt-1
		G(n+1,0) = Coarse(G(i,0))
		U(n+1,0) = G(n+1,0)
	endfor
	
	# Parareal iteration until K
	
	for k=1:K
	
		# execute fine method in parallel
		
		for n=k-1:Nt-1	
			F(n+1,k-1) = Fine(U(n,k-1))
			Delta = F(n+1,k-1) - G(n+1,k-1)
		endfor
		
		# Parareal update procedure in serial
		
		for n=k-1:Nt-1
			G(n+1,k) = Coarse(U(n,k))	
			U(n+1,k) = G(n+1,k) + Delta
		endfor
	
	endfor
\end{lstlisting}

FESOM2 is a called by a command line operation, for which the location, executable and available resources have to be defined. The evolution of each initial value problem, in serial and parallel, within the sbash script is carried out by the \texttt{srun} command: \\
\lstset{basicstyle=\fontsize{9}{13}\selectfont\ttfamily}
\begin{lstlisting}
	srun $flgs --chdir=$FOLDER_PATH ./fesom.x > "$FOLDER_PATH/fesom2.0.out" &
\end{lstlisting}
\vspace{.4cm}
Here, the value  \texttt{\$flgs} contains the amount of nodes and processors requested for the respective run and \texttt{\$FOLDER\_PATH} denotes the location of the folder on the cluster in which FESOM2 is executed. For each coarse, fine and iterative solution folders were created to store simulation data. The ampersand at the end of the \texttt{srun} command is set to run the sub-job in parallel with other instructions. For the serial evolution of the coarse propagator the ampersand is omitted. 

For the evaluation of the correction procedure by the fine and coarse propagator at a given Parareal iteration $k$, the operator collection \textit{NCO} was used, see \hyperlink{nco}{http://nco.sourceforge.net/}. It allows for manipulation of model data by command line tools and is designed to handle netCDF files. During the execution of the Parareal algorithm, one has to compute the difference and sum of the variables stored in two netCDF files. During the time-parallel computation of the fine propagator, the difference to the coarse propagator at each time slice is calculated, see line 17 in Lst. \ref{lst:Parareal}. In the sbatch script, the NCO operator \texttt{ncbo} is used as: \\
\begin{lstlisting}
	ncbo -O --op_typ=diff $PATH_FINE_FILE $PATH_COARSE_FILE $PATH_FILE_OUT
\end{lstlisting}
\vspace{.4cm}
The \texttt{ncbo} operator is given the type of operation by \texttt{--op\_typ=} as well as input and output file locations. The output file, representing the difference term  \texttt{Delta} in the pseudo-code, is stored separately and will be used in the sum operation during the algorithm's update step. Accordingly, the addition required in the update procedure (line 24 in Lst. \ref{lst:Parareal})  is represented by: \\
\begin{lstlisting}
	ncbo -O --op_typ=add $PATH_COARSE_FILE $PATH_DELTA $PATH_ITER_FILE
\end{lstlisting}
\vspace{.4cm}
With these command line instructions within the sbatch script, the core tasks of Parareal can be submitted to the DKRZ cluster. The script additionally contains commands to store simulation data of the Parareal approximation after each iteration for post-processing vizualisation and error estimates.

\section{Numerical Experiments}

In this section settings for the numerical experiments are introduced. The overall goal of this study was to investigate appropriate numbers of temporal sub-intervals that allow for parallelization in time by the classical Parareal algorithm. Before we dive into the details of the numerical experiments, spatial parallelization and runtime estimates for serial computations with different time step-sizes are presented. 

\subsection{Runtime Estimates}

Investigation of the spatial parallelization capabilities of FESOM2 allows for an optimal choice of processors to run FESOM2 on the PI mesh. For runtime evaluations with respect to time stepping, we are restricted, by the design of FESOM2, to the time integration method Adams-Bashforth of third order. This restricts the distinction of fine and coarse propagators to the variation of the time step-sizes. 

\begin{figure}
	\centering
	\subfloat{\includegraphics{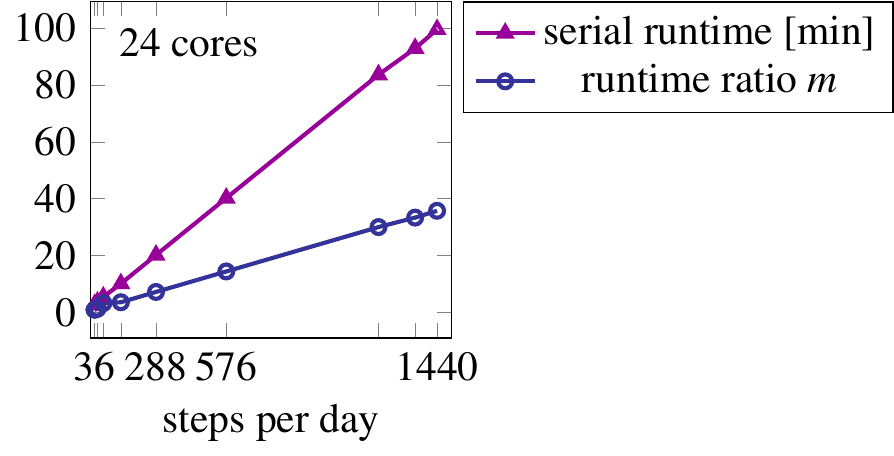}}
	\subfloat{\includegraphics{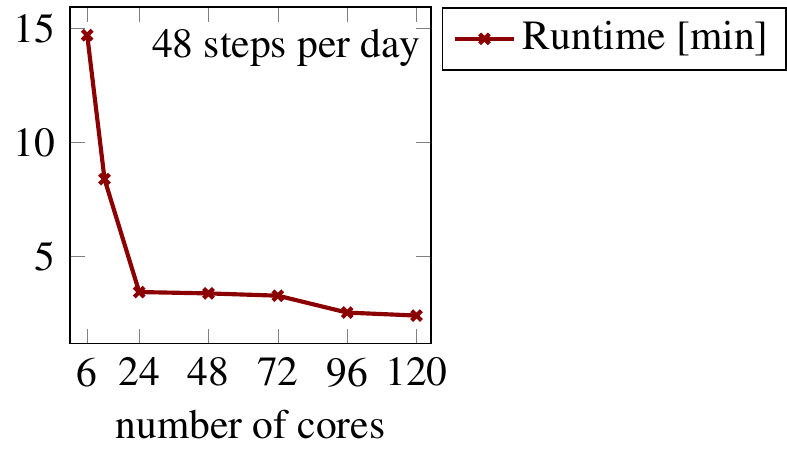}}
	\caption{Left: Runtime evaluations on the PI mesh computed over $T=1$ year simulation time for different time steps sizes on 24 processors. Reference for the runtime ratio $m$ is a serial run with 36 steps per day. Right: Runtime evaluation with fixed time step-size and for various amount of parallel processes.}
	\label{FIG:FESOM2scale}
\end{figure}

For each time slice,  the methods will be executed with the time step-sizes $\Delta_G$ for the coarse and $\Delta t_F$ for the fine propagator. In Fig. \ref{FIG:FESOM2scale}, the runtime of FESOM2 with respect to the used time step-size for a fixed amount of processors are depicted. Due to the low spatial resolution of the PI mesh, the threshold for spatial domain decomposition in serial runs is found at 24 processors, as it is shown in the right picture. Thus, for each propagator run in Parareal, we applied 24 cores for spatial domain decomposition to ensure efficient execution of FESOM2. This evaluation 
ensures that the simulations were carried out with acceptable space-parallel efficiency. For performance evaluations of the code on finer resolved meshes the reader is referred to \cite{FESOM2performance}. 

In FESOM2, the time step-size is defined by steps per day (\textit{spd}), where 36 spd corresponds to 40 minutes or 2400 seconds for one time-step. The number of steps per day for the coarse and fine method are denoted as $N_G$ and $N_F$, respectively. Running FESOM2 with 36 spd executes the model at the limit of the CFL condition on the PI mesh, as we found during first numerical tests. Therefore, this value represents a threshold for a time step increase. 
Considering the given speed-up estimate in Eq. \eqref{EQ:speedup}, we have chosen the coarse propagator to proceed with the allowed minimum of steps per day, while the fine propagator was set to integrate in time as fine as possible. In the left image of Fig. \ref{FIG:FESOM2scale}, the runtime ratio $m = \tau_F / \tau_{N_G = 36}$ is shown with respect to the fastest setting with 36 steps per day. The execution time ratio should be chosen as $m \gg 2$ to be able to generate speed-ups according to the estimate in Eq. \eqref{EQ:speedup}.

We decided to limit the number of steps by pragmatic considerations. A time step reduction to 1440 spd increases the runtime by an estimated factor of 40, which can be considered as a decent runtime ratio $m$. From the scalability test results depicted in Fig.\ref{FIG:FESOM2scale}, we can confirm this holds true. Any further time step refinement did not appear useful to us and therefore remained unconsidered. 

Ultimately, the finest allowable time step-size in FESOM2 denotes 1 second, corresponding to 86400 spd, which is not computationally feasible anymore, as it increases the runtime ratio to $m \approx 2400$. We remark that FESOM2 requires  time step-sizes  to be an integer divisor of 86400, which is equal to the amount of seconds within one day.  

\subsection{Error Estimates}

To estimate the deviation in solutions due to different time step-sizes, two error norms are introduced: the euclidean and maximum error norms. The errors are evaluated according to the respective serial fine run solution, here denoted by  $u^F_{N_t}$, with $N_F$ as:
\begin{equation}
	\begin{split}
		E^{N_F}_\infty \; &= \; \frac{\| U^k_{N_t} - u^{N_F}_{N_t} \|_\infty}{\| u^{N_F}_{N_t} \|_\infty},  \\
		E^{N_F}_2 \; &= \; \frac{\| U^k_{N_t} - u^{N_F}_{N_t} \|_2}{\| u^{N_F}_{N_t} \|_2}.
	\end{split}
	\label{EQ:metrics}
\end{equation}
Since all experiments were carried out with 36 spd for the coarse solver, the computed norms are given with respect to the fine method's time step-size, denoted by the superscript $N_F$. The subscript $N_t$ indicates that the error norms were evaluated at the endpoint of the simulation time $t_{N_t} = T$. 

All error estimates considered solely the surface layer of the ocean. The dynamics of the so-called deep ocean are reduced in comparison to the upper layers. We found for all test cases that the Parareal algorithm showed better convergence behavior with increasing ocean (layer) depth. Accordingly, we chose the ocean surface layer as an appropriate indicator for success of the time-parallel approach. 

Both norms provide different information to the approximation process of Parareal. The euclidean norm estimates the overall deviation of two fields to each other, while the maximum norm corresponds to local errors. With observation of both we aim for a more holistic understanding and evaluation of the Parareal algorithm applied to FESOM2.

\subsection{Numerical Test Cases}

To test the convergence behavior of Parareal, the algorithm was carried out for three test cases of FESOM2. We distinguished the experiments by terms of time slice length $\Delta T$. For all test cases the coarse propagator was executed with 36 spd corresponding to a 2400 seconds time step-size $\Delta t_G$. Consequently, the smallest time slice length $\Delta T$ available  is $\Delta T = \Delta t_G$. For the  first numerical experiment in this study we started with the smallest time slice length. To estimate the convergence behaviour we increased $\Delta T$ to one day and month for the second and third case, respectively. The simulation time interval $T$ was split into $N_t=12$ time slices of equal length in order to approximate one year of climate simulation in the third test. With keeping $N_t$ fixed the different $\Delta T \in \{ 2400 \, \text{seconds} , 1 \, \text{day} , 1 \, \text{month} \}$ result in the time intervals $T \in \{ 8 \, \text{hours}, 12 \, \text{days} , 1 \, \text{year} \}$. The time step size $\Delta t_F$ for the fine method was refined by multiples of 2, such that the number of steps per day increased as $N_F = 72, 144, 288$. This corresponds to an estimated runtime ratio of $m = 2,4,8 $. The step-size refinement for the fine method by multiples of 2 was determined by the design of FESOM2. Permissible number of steps per day are integer divisors of 86400. Accordingly, we organized the numerical test settings into \\

\begin{table}[H]
	\centering
	\caption{Numerical experiment settings for Parareal.}
	\begin{tabular}{|l|c|c|c|c|c|}
		\hline
		Experiment No. & $\Delta T$ & $T$ & $N_F$ & $N_G$ & $N_t$ \\
		\hline
		1 & 2400 seconds & 8 hours & {} & {} & {} \\
		2 & 1 day & 12 days & $72, 144, 288 $ & 36 &12\\
		3 & 1 month & 1 year & {} & {} & {} \\
		\hline
	\end{tabular}
\end{table}

\subsection{Results}

In this section the results of the three numerical experiments are presented. The stability of Parareal when applied to FESOM2 is dependent on $N_F$, although there is no generalized statement possible. For evaluations of the time-parallel algorithm with a time slice length $\Delta T = \Delta_G = 2400$ seconds we found the limit for a stable iteration at 288 spd for the fine method. In the second experiment a time slice length of one day would allow for fine time step-sizes up to 1440 spd. For the last test case $\Delta T = 1$ month the Parareal algorithm suffered from blow-ups for settings with 144 spd and higher. 
During this study we considered two perspectives of a successful implementation. From a pragmatic point of view, we defined the error threshold $\varepsilon \leq 10^{-2}$ for both norms as successful. However, in terms of a thorough analysis and to fully assess convergence behavior we applied Parareal for $K=N_t$ in order to check for convergence to machine precision. 

\subsubsection{Experiment 1}

We will begin with introducing the results for the shortest time slice length $\Delta T = \Delta_G = 2400$ seconds, corresponding to one time step over $\Delta T$ for the coarse propagator.

\begin{figure}[H]
	\centering
	\subfloat{\includegraphics{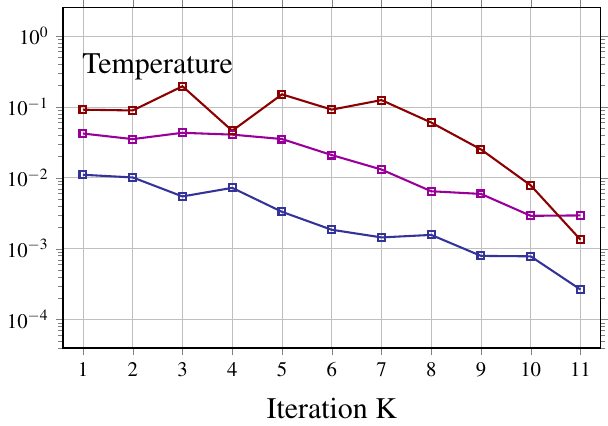}
	}\subfloat{\includegraphics{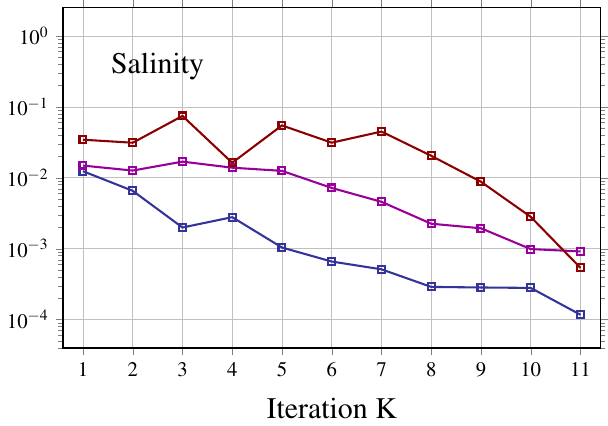}}

	\subfloat{\includegraphics{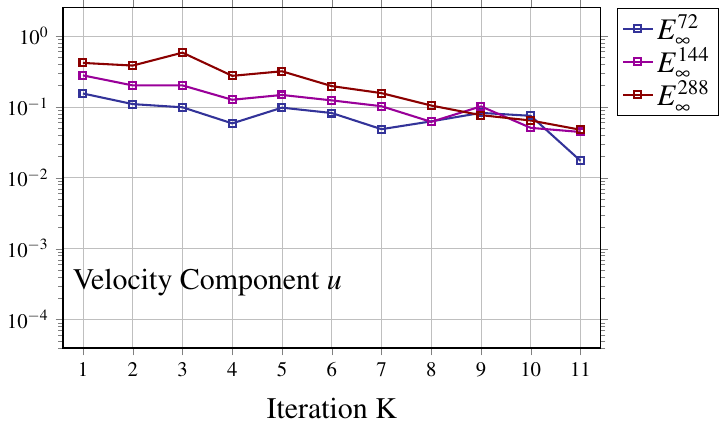}}	
	\caption{Parareal convergence results for $T = 8$h with $\Delta T = \Delta_G = 2400$ seconds ($N_C = 36$ spd) and $N_F =72, 144, 288$. The relative maximum norm $E^{N_F}_\infty$ was evaluated at the final simulation time and are depicted in dependency of the Parareal iteration count $K$.}
	\label{FIG:p36_sec_inf}
\end{figure}

\begin{figure}
	\centering
	\subfloat{\includegraphics{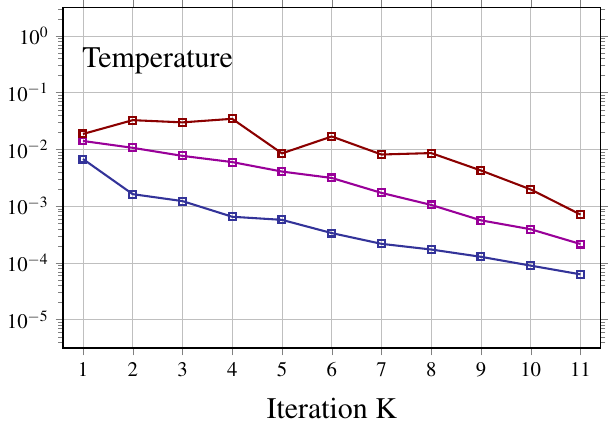}
	}\subfloat{\includegraphics{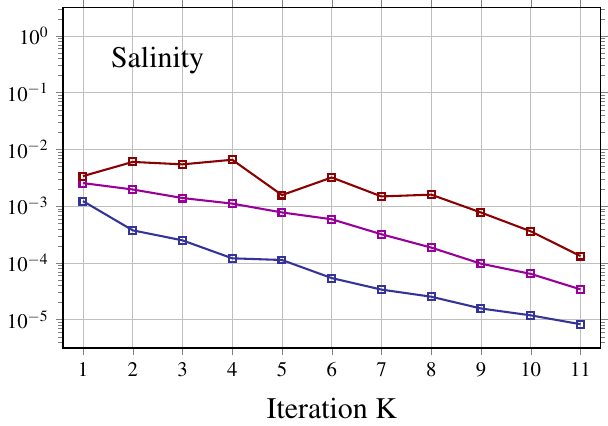}}

	\subfloat{\includegraphics{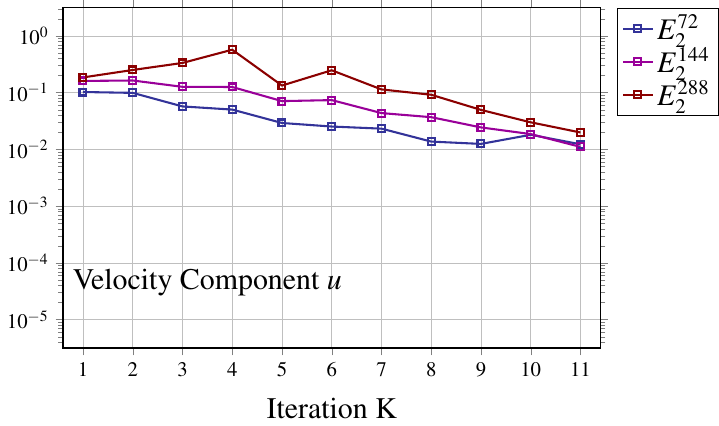}}
	\caption{Parareal convergence results for $T = 8$h with $\Delta T = \Delta_G = 2400$ seconds ($N_C = 36$ spd) and $N_F =72, 144, 288$. The relative euclidian norm $E^{N_F}_2$ was evaluated at the final simulation time and are depicted in dependency of the Parareal iteration count $K$.} 
	\label{FIG:p36_sec_eucl}
\end{figure}

In the Figs.\ref{FIG:p36_sec_inf} and \ref{FIG:p36_sec_eucl} the relative maximum and euclidean norm are depicted for the variables temperature, salinity and velocity component $u$. For the error estimation serial fine computations were conducted with respect to choices of the fine method during the respective Parareal run. Both metrics were evaluated for each iteration $k$ at the final simulation time.
From Fig.\ref{FIG:p36_sec_inf} a distinctive separation in the error curves with respect to the fine propagator is observed for temperature and salinity. Parareal convergence to floating point precision could not be achieved in any of the configurations. An error reduction of two orders for temperature and salinity is obtained only by the last two iterations. The results for the zonal velocity $u$ show a less distinctive separation with respect to the choice of the fine method. With increasing iteration count $k$ error curves seem to align. An error reduction for $u$ of order one requires 11 Parareal iterations for all three settings of $N_F$. For the variables governed by the momentum equations the algorithms performance is evidently stagnating. The error norm of the last iteration $k=12$ is omitted in all figures since the algorithm recomputes the fine serial solution, for which the error estimate self-evidently is zero. Refinement of $N_F$ further decreases the convergence behavior until the algorithm becomes unstable for $N_F > 288$ spd. \\

From Eq.\eqref{EQ:speedup} we recall the speedup bounds for the three sub-settings $m =2, 4, 8$, resulting in $S \leq\min \left\{ \frac{m}{k+1} , \frac{12}{k} \right\}$. For $m=2$ we find the speedup bound is already limited to $S \leq 1$ by the first Parareal iteration $k=1$. Recall that the speedup bound neglects all communication cost on the HPC. Hence, the test runs with $N_F=72$ spd were conducted for the sake of completeness of the convergence examination. With $N_t=12$ the speed up limit is governed by the runtime ratio $m$, yielding limits for allowable maximal iteration count by $K \leq0, 2, 6$ for $N_F =72, 144, 288$, respectively. Even under consideration of the forgiving error threshold of $\epsilon \leq 10^{-2}$ we observed the required Parareal iterations to exceed the given limit.
Relative euclidian error norm are given in Fig.\ref{FIG:p36_sec_eucl} for which the convergence behavior is comparable to the results of the relative maximum norm in Fig.\ref{FIG:p36_sec_inf}. For temperature and salinity the euclidian error is found to be shifted by a magnitude of one, whereas the outcome for the velocity component $u$ shows almost no shift. Again, convergence to floating point precision was not achieved. For all iterations the error curves for salinity are below $10^{-2}$, suggesting that the Parareal approximation and serial computation salinity fields are in acceptable agreement. Convergence in terms of the threshold $\epsilon \leq 10^{-2}$ was achieved after one iteration for salinity. The temperature fields reached the threshold after 1,3 and 7 iterations for $N_F=72, 144, 288$, respectively. Resulting in theoretical speedups of $S \leq 1$ for all fine propagator configurations. The convergence for the zonal velocity did stagnate and we were not able to reduce the error significantly during the time-parallel iteration procedure. 

\subsubsection{Experiment 2}

In Figs.\ref{FIG:p36_day_inf} and \ref{FIG:p36_day_eucl} the results are shown for $\Delta T = 1$ day, resulting in a simulation time of $T=12$ days. For the relative maximum norm depicted in Fig.\ref{FIG:p36_day_inf} we observe the stagnation of Parareal convergence for temperature and salinity, while the iteration for the zonal velocity $u$ diverges for the first iterations before a a minor reduction of error could be achieved. Destinction of the error curves with respect to the fine propagator time step is only observable for the divergent behavior in the zonal velocity. The pre-defined threshold of $10^{-2}$ was not reached except in the last iterations of Parareal for salinity and temperature. 
We observed that the relative euclidean norm for temperature and salinity was below the threshold for all iterations. At the same time no significant error reduction was achieved and the error estimates for $u$ show a divergent behavior similar to the relative maximum norm. We found the performance of Parareal to further impair with increasing the time slice length $\Delta T$ and refinement of the fine propagator. 

The initial relative euclidean error for temperature and salinity appears to be acceptable without applying Parareal, while the error of the velocity fields does not allow for any speedup generation by the algorithm. Although, the governing equations Eq.\eqref{EQ:tracer} for temperature and salinity are dependent on the velocity, the convergence seem to be less affected by the fine methods time step-size $\Delta t_F$ and the time slice length $\Delta T$. 

\afterpage{%
	\newgeometry{left=3cm,right=3cm,top=1cm,bottom=1cm,includeheadfoot}
	% material for this page

\begin{figure}
	\centering
	\subfloat{\includegraphics{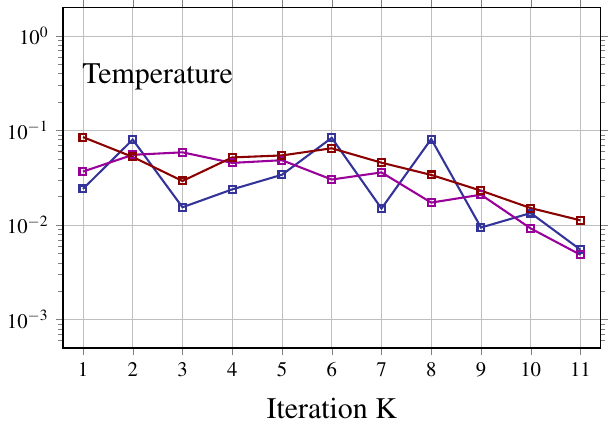}
	}\subfloat{\includegraphics{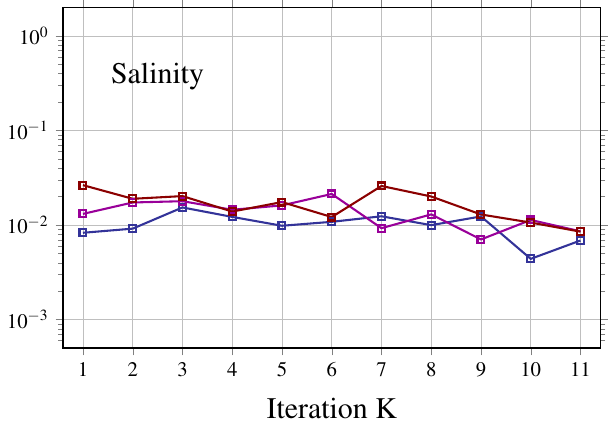}}

	\subfloat{\includegraphics{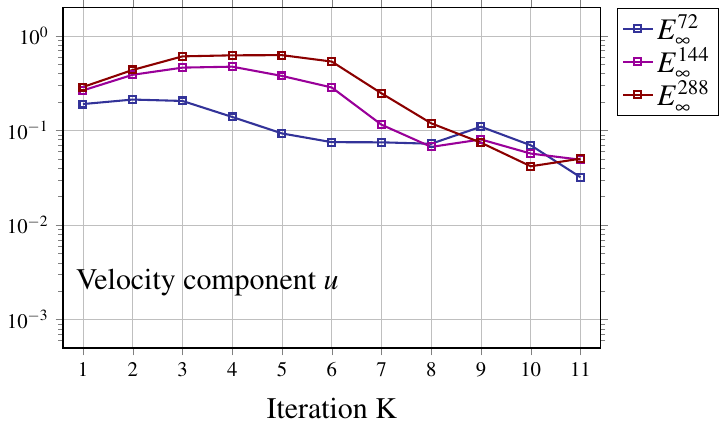}}
	\caption{Parareal convergence results for $T = 12$ days, $\Delta T = 1$ day, $N_C = 36$ spd, $N_F \in [72,144,288]$ spd. The relative maximum norm $E^{N_F}_\infty$ is depicted in dependency of the Parareal iteration count $K$.}
	\label{FIG:p36_day_inf}
\end{figure}

\begin{figure}[H]
	\centering
	\subfloat{\includegraphics{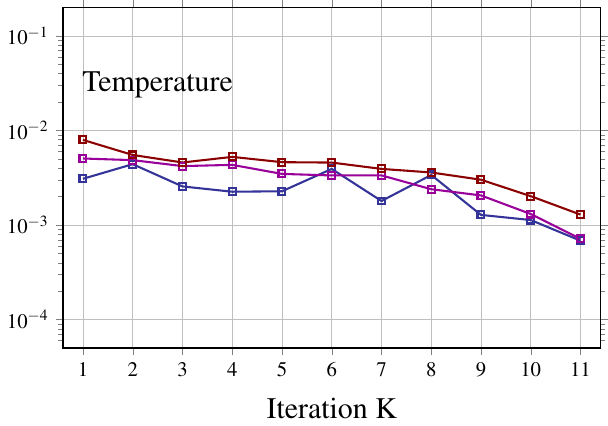}
	}\subfloat{\includegraphics{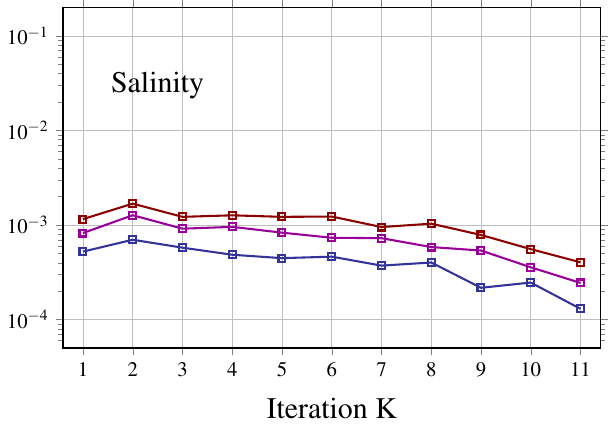}}

	\subfloat{\includegraphics{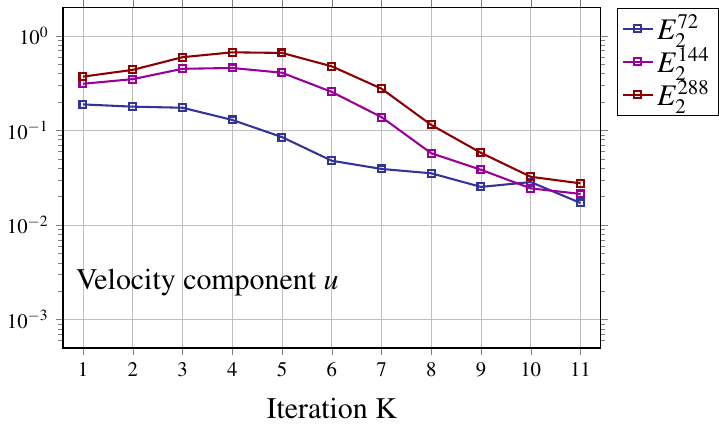}}
	\caption{Parareal convergence results for $T = 12$ days, $\Delta T = 1$ day, $N_C = 36$ spd, $N_F =72,144,288$ spd. The relative euclidian norm $E^{N_F}_2$ is depicted in dependency of the Parareal iteration count $K$.}
	\label{FIG:p36_day_eucl}
\end{figure}

	\restoregeometry\clearpage	
}
\newgeometry{left=3cm,right=3cm,top=2cm,bottom=2cm,includeheadfoot}

\subsubsection{Experiment 3}

\begin{figure}[H]
	\centering
	\subfloat{\includegraphics{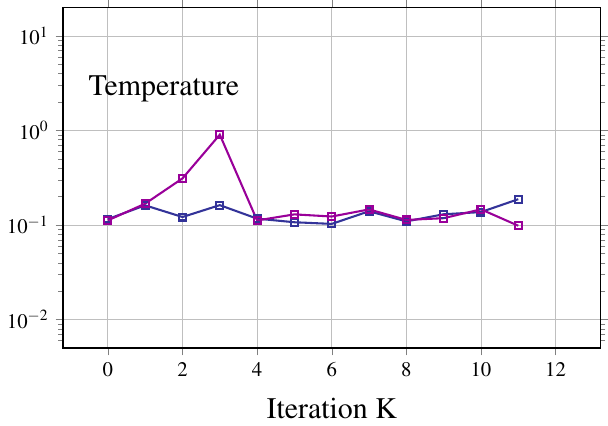}
	}\subfloat{\includegraphics{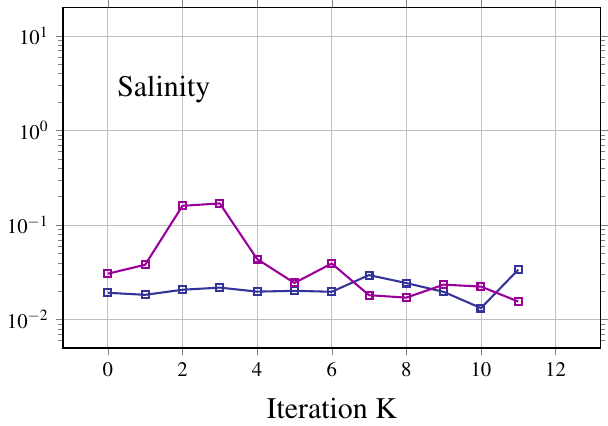}}

	\subfloat{\includegraphics{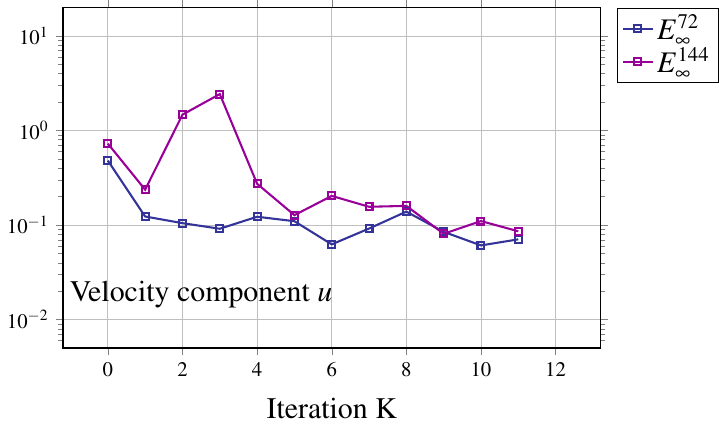}}
	\caption{Parareal convergence results for $T = 12$ months with $\Delta T = 1$ month, $N_C = 36$ spd and $N_F =72,144$ spd. The relative maximum norm $E^{N_F}_\infty$ was evaluated at the final simulation time and are depicted in dependency of the Parareal iteration count $K$.}
	\label{FIG:p36_mon_inf}
\end{figure}

\begin{figure}[H]
	\centering
	\subfloat{\includegraphics{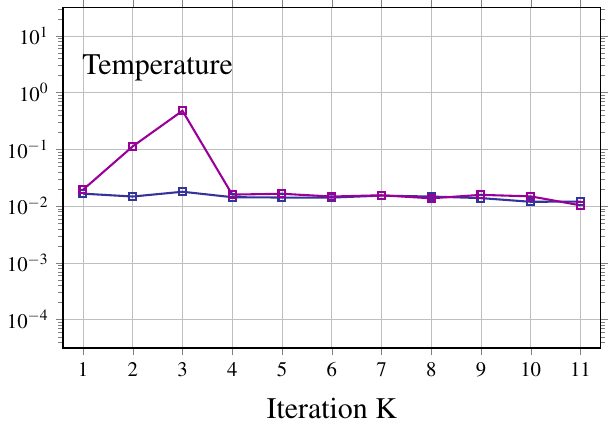}
	}\subfloat{\includegraphics{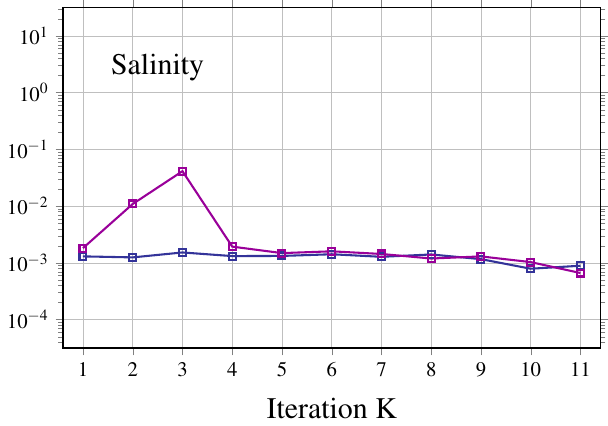}}

	\subfloat{\includegraphics{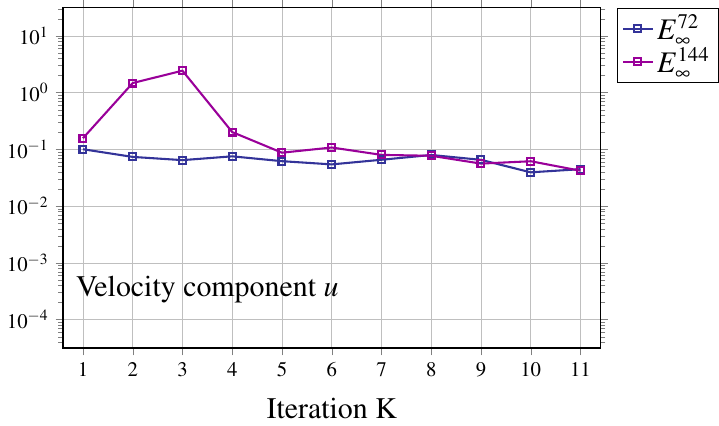}}
	\caption{Parareal convergence results for $T = 12$ months with $\Delta T = 1$ month, $\Delta_G = 36$ spd and $\Delta t_F \in 72,144$ spd. The relative euclidian norm $E^{N_F}_2$ was evaluated at the final simulation time and are depicted in dependency of the Parareal iteration count $K$.}
	\label{FIG:p36_mon_eucl}
\end{figure}

For the last test case we applied the Parareal algorithm to a simulation time of $T= 1$ year, split by $N_t = 12$ time slices of $\Delta T= 1$ month each. The evaluation of both error norms is given in Figs.\ref{FIG:p36_mon_inf} and \ref{FIG:p36_mon_eucl}. It has emerged that Parareal was not able to provide a stable iteration for time step-sizes $N_F > 72$. By the speedup estimate in Eq.\eqref{EQ:speedup} we know that for the configuration with $N_F = 72$ a runtime reduction is not obtainable in the first place. 

While the error residual stagnates for both error metrics in case of $N_F = 72$, we find the Parareal iteration for $N_F = 144$ to suffer from a blow-up at iteration $k=3$. A blow-up in the Parareal iteration will force the FESOM2 solver to terminate its computation. If the CFL condition is violated due to diverged initial states provided by Parareal, the simulation of a initial value problem over a given time slice cannot be carried out. Eventually, the necessary subtraction and summation operations cannot be executed since the required input data files were not created. With the Parareal instructions given in a sbatch file these tasks were omitted and a runtime error notification is written to a log file. Nevertheless, the script is able to continue since the file containing $G(U^{k+1}_n)$ is modified by adding the correction. In case of no correction the unedited coarse solution file is used as $U^{k+1}_{n+1}$. Although the script continues no error reductions can be achieved in any of the settings. Consequently, $N_F=288$ was omitted since no further insights could be provided.
The shift of one magnitude in error between the relative maximum und eculidean metrics suggests that most parts of the variable fields are acceptable in error only under consideration of the threshold. Parareal appeared incapable of reducing the differences in some parts of the fields, where the impact of time step variation appears to be significant. 
With further increased time slice length $\Delta T$ and refinement of the temporal resolution for the fine method we observed an additional decrease in the performance of Parareal in all test cases. Unfortunately, even for the shortest time slice length $\Delta T = \Delta_G$, where the coarse propagator executes one step in time, the algorithm fails to approximate a solution to the serial fine run. Convergence to double precision is given only if Parareal executes all iterations $K=N_t$ and thereby reproduces the respective serial simulation. Under consideration of the quite tolerant threshold $\varepsilon$, convergence by the time-parallel algorithm is prevented by the horizontal velocity components governed by the momentum equation for all test cases. Both error norms in the velocity fields could not be significantly reduced by Parareal during the experiments. 

\newpage

\section{Discussion}

For the assessment of the presented performance results we orientated the discussion along the following points:
\vspace{.4cm}
\begin{enumerate}
	\item The model equations' sensitivity to initial values.
	\item The technical design of FESOM2 when it comes to storing results and the continuation of simulation runs.
	\item FESOM2 ability to simulate long-term ocean dynamics with respect to averaged diagnostic variables.
\end{enumerate}
\vspace{.4cm}

Firstly: the simulation of the ocean dynamics on the discretized globe needs support by models, e.g. for turbulence, in order to provide a stable solver that allows for appropiate large time step-sizes. Although, the complexity is somewhat reduced with respect to fully resolved turbulent flows, one has to take the Navier-Stokes equations' sensitivity to initial values into account. The tracer equations for temperature and salinity are not prone to initial sensitivity as the momentum equations, but are dependent on the velocity components. Accordingly, we observed better convergence behaviour of temperature and salinity in comparison to the velocities. Although improved, even an convergence order of 1 could not be obtained for the tracer variables. 
Furthermore, with the increase of time slice length $\Delta T$ throughout the numerical experiments we observed an error build-up that surpasses the ability of Parareal to reduce error by iteration, especially for the horizontal velocity. Unfortunately, the stagnant convergence behavior of the horizontal velocity during the iteration procedure of Parareal prohibited speedup generation for all test cases. E.g., the configuration with $\Delta T = 1$ day and $N_F = 144$ the error norms for temperature and salinity fell below the threshold $\varepsilon$ within a reasonable amount of iterations. 

\begin{figure}[H]
	\centering
	\subfloat{\includegraphics{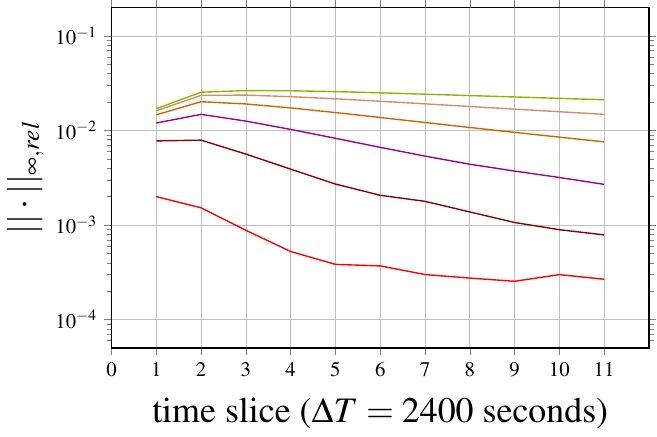}
	}\subfloat{\includegraphics{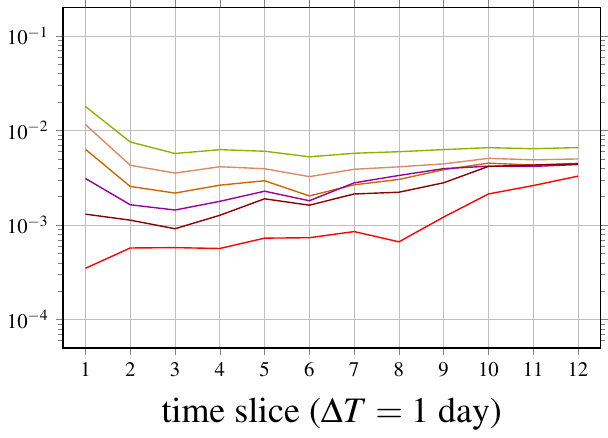}}

	\subfloat{\includegraphics{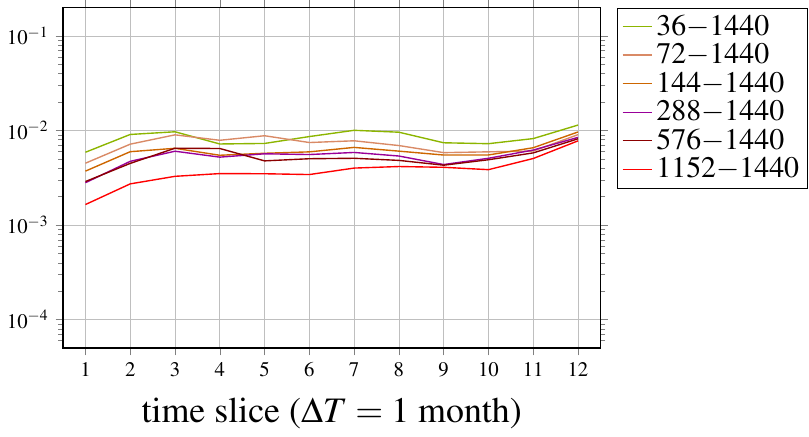}}
	\caption{Evaluation of the relative maximum norm for serial runs in $36,72,144,288,576,1152$ spd with respect to the reference solution with $1440$ spd. The norm was evaluated at each time slice of varying length for the time-averaged temperature for each time slice $\Delta T$.}
	\label{FIG:FESOM2serialnorm}
\end{figure}

Secondly: FESOM2 stores prognostic variables at the latest point in time in restart files with successful termination of a simulation run. From the restart files one can continue to run the model. Unfortunately, FESOM2 introduces small deviations to the initial time stepping process. Which means, that a consecutive execution over 2 years of simulation time results in different solutions than a simulation of same length that is restarted after one year. For the sake of comparison we had to carry out the serial reference runs for Parareal in a restarted manner. Which shows that FESOM2 is not meant to preserve numerical accuracy when performing restarts. Therefore, convergence to a consecutively carried out simulation was not possible in the first place, as we discovered during our evaluation process. Which is leading us directly to the last point:

One needs to recall that FESOM2 aims to model the large-scale ocean dynamics and sea-ice evolution over years up to decades. These models are developed and tuned to approximate the physical phenomenons according to the spatial resolution of the respective mesh. The choice of time step-sizes appears to be a possibility to ensure numerical stability in simulations with FESOM2 rather than to achieve numerical accuracy in time. In Fig.\ref{FIG:FESOM2serialnorm} we compared serial runs of different amounts of steps per day to a reference solution with 1440 spd. All runs were executed in a restarted manner after each time slice. The time slices were chosen in accordance with the runtime evaluation setup, depicted in Fig.\ref{FIG:FESOM2scale}.

The finest temporal resolution of $1440$ spd was used as the reference solution to compute the relative maximum norm at all time slices and at the surface layer of the ocean. The variables were averaged in time over each time slice. We restricted this investigation to errors in temperature, since the depicted results are representative in their tendency for salinity and horizontal velocity. 

With increasing the time slice length $\Delta T$ a decline in the significance of the chosen time step-sizes was observed. Instead the errors to the reference solution for all runs appear to align around $10^{-2}$ with increasing simulation time. It was not clear to us, whether this behavior is due to the atmospheric forcing and climatology data sets used for the simulations or how the model parameter tuning impacts the results. Consequently, we failed being able to make a statement about the approximation accuracy with respect to the time step-size. With longer simulation times, for which FESOM2 is intended for, all runs appear to lie within the same error range with respect to the reference solution. 

\section{Conclusion}

This paper investigates the impact of time step-sizes and simulation intervals on the Parareal algorithm applied to the ocean-circulation and sea-ice model FESOM2. Besides an experimental convergence evaluation this study was conducted to assess the implementation to a state-of-the-art climate model. 
We summarize the findings of the study by: \\

\begin{itemize}
	\setlength\itemsep{1em}
	\item Parareal convergence to floating point precision was obtained only with iteration $K=N_t$, where the serial reference solution is re-computed.
	\item Based on the comparison of time averaged variables the desired configuration of FESOM2 with respect to time integration is the largest stable time step-size available. Smaller time step sizes do not contribute to a better prediction of diagnostic variables. 
	\item Hence, reducing time step sizes in order to define a fine method for Parareal only artificially increases wall clock times of FESOM2.

\end{itemize} 
\vspace{.4cm}

In the here presented numerical study of the Parareal algorithm applied to FESOM2 we were not able to generate speedups for variety of time slice lengths and fine propagator configurations. Convergence to floating point precision was only achieved when the maximum iteration count $K=N_t$ was reached and the respective fine serial solution was recomputed. Even for the proposed pragmatical error threshold $\varepsilon \leq 10^{-2}$ the algorithm would exceed the limiting iteration count, that would allow for at least minimal runtime reduction. We conclude that the application of Parareal in order to generate speedups is not applicable for two main reasons: The algorithms struggles to converge fast enough for runtime reduction and secondly, the informational value by using a more accurate fine propagator in time is neglectable. Variations in the local states in time do occur but play no significant role in the time-averaged quantities of the long-term ocean dynamics. 

In the authors opinion the concept of Parareal can still be useful. The approach of using meshes of different spatial resolution for the coarse and fine propagator may be an option. For future work one would have to face the main challenges of spatial coarsening for FESOM2 meshes, which lies in the generation of a suitable coarser mesh and for interpolation methods conserving mass and momentum. Further, one will likely have to optimize the model parameters for the coarse mesh setting in order to provide comparable ocean dynamics. The overall goal for a micro-macro variant of Parareal with two spatial resolutions will be the approximation of time-averaged diagnostic variables instead of converging to solutions to specific points in time. 
There exist successful implementations of the so-called micro-macro Parareal algorithm for an energy-balance-model \cite{Slawig2018} and the time-parallel simulation of homogeneous turbulence \cite{Lunet2018}. Especially the former study provides the property to preserve the fine scale contributions of the fine resolution during the iteration procedure. In combination with the design of FESOM2 to approximate ocean dynamiscs based on the spatial resolution we assume a successful application of the micro-macro Parareal algorithm to FESOM2 possible.

%%% SOURCES %%%

\end{document}